\documentclass[reqno,12]{amsart}
\usepackage{amssymb}
\usepackage{amsmath, amsthm}
\usepackage{hyperref}
\usepackage{enumerate}
\usepackage{graphicx}
\usepackage{bm}
\newcommand{\h}{\frac{1}{2}}
\AtBeginDocument{\noindent\small
\vspace{9mm}}
\begin{document}
\title[On an identity for H-function]{On an identity for H-function }
\author[Arjun K. Rathie ]{Arjun K. Rathie}
\address{Department of Mathematics, School of Mathematics and Physical sciences, Central University of Kerala, Periye P.O., Kasaragod- 671316, Kerala, INDIA}
\email{akrathie@cukerala.ac.in }
%\address{}
%\thanks{Submitted:  2017}
\begin{abstract}
The main objective of this research note is to provide an identity for the H-function, which generalizes two identities involving H-function obtained earlier by Rathie and Rathie et al.
\end{abstract}
\maketitle
\numberwithin{equation}{section}
\newtheorem{theorem}{Theorem}[section]
\newtheorem{lemma}[theorem]{Lemma}
\newtheorem{proposition}[theorem]{Proposition}
\newtheorem{corollary}[theorem]{Corollary}
\newtheorem{Definition}[theorem]{Definition}
\newtheorem*{remark}{Remark}
\newtheorem{integral}[theorem]{Integral}
\section{Introduction}
In 1981, Rathie [2] established the following identity for the H-function viz.
\begin{align}
H_{\;p+1,\; q+1}^{\;m,\; n}& \left[ z |\begin{array}{c} \;{}_1(a_j, e_j)_p, (\alpha, \lambda)\\ 
\;{}_1(b_j, f_j)_q, (\alpha, \lambda)\end{array} \right] \nonumber\\
&  = \frac{1}{2\pi i}\left\{ e^{i \pi \alpha} H_{\;p,\; q}^{\;m,\; n} \left[ ze^{-i\pi \lambda} |\begin{array}{c} \;{}_1(a_j, e_j)_p\\  \;{}_1(b_j, f_j)_q\end{array} \right] \right. \nonumber\\
& \qquad \qquad - \left. e^{-i \pi \alpha} H_{\;p,\; q}^{\;m,\; n}  \left[ ze^{i\pi \lambda} |\begin{array}{c} \;{}_1(a_j, e_j)_p\\  \;{}_1(b_j, f_j)_q\end{array} \right] \right\}
\end{align}
Very recently, Rathie et al.[3] established another identity for the H-function viz.
\begin{align}
{H}_{\;p+1,\; q+1}^{\;m+1,\; n+1}& \left[ z |\begin{array}{c} (\alpha, \lambda),\;{}_1(a_j, e_j)_p\\ 
 (\alpha, \lambda), \;{}_1(b_j, f_j)_q\end{array} \right] \nonumber\\
&  =  e^{i \pi \alpha}\; H_{\;p+1,\; q+1}^{\;m+1,\; n+1} \left[ ze^{-i\pi \lambda} |\begin{array}{c}  (2\alpha, 2\lambda),\;{}_1(a_j, e_j)_p\\ (2\alpha, 2\lambda), \;{}_1(b_j, f_j)_q \end{array} \right]  \nonumber \\
& \qquad +  e^{-i \pi \alpha}\; H_{\;p+1,\; q+1}^{\;m+1,\; n+1}  \left[ ze^{i\pi \lambda} |\begin{array}{c}(2\alpha, 2\lambda), \;{}_1(a_j, e_j)_p\\ (2\alpha, 2\lambda), \;{}_1(b_j, f_j)_q \end{array} \right] 
\end{align}
For interesting applications of the identities (1.1) and (1.2), we refer the recent paper by Rathie et al.[3]\\
The aim of this short note is to provide a natural generalization of (1.1) and (1.2).
\section{Main Result}
The identity for the H-function to be established in this note is  the following. 
\begin{align}
2\pi i \; &  H_{\;p+2,\; q+2}^{\;m+1,\; n+1} \left[ z |\begin{array}{c} (\beta, \delta),\;{}_1(a_j, e_j)_p, (\alpha, \lambda)\\ (\beta, \delta), \;{}_1(b_j, f_j)_q, (\alpha, \lambda)\end{array} \right] \nonumber\\
&  =  e^{i \pi (\alpha+\beta)}\; H_{\;p+1, \; q+1}^{\;m+1,\; n+1} \left[ ze^{-i\pi (\lambda+\delta)} |\begin{array}{c} (2\beta, 2\delta), \;{}_1(a_j, e_j)_p\\ (2\beta, 2\delta), \;{}_1(b_j, f_j)_q\end{array} \right]  \nonumber\\
&  +  e^{i \pi (\alpha-\beta)}\; H_{\;p+1,\; q+1}^{\;m+1,\; n+1}  \left[ ze^{-i\pi (\lambda-\delta)} |\begin{array}{c} (2\beta, 2\delta), \;{}_1(a_j, e_j)_p\\ (2\beta, 2\delta), \;{}_1(b_j, f_j)_q\end{array} \right] \nonumber\\
& -  e^{-i \pi (\alpha-\beta)}\; H_{\;p+1,\; q+1}^{\;m+1,\; n+1}  \left[ ze^{i\pi (\lambda-\delta)} |\begin{array}{c} (2\beta, 2\delta), \;{}_1(a_j, e_j)_p\\ (2\beta, 2\delta), \;{}_1(b_j, f_j)_q\end{array} \right]\nonumber\\
&  - e^{-i \pi (\alpha+\beta)}\; H_{\;p+1,\; q+1}^{\;m+1,\; n+1} \left[ ze^{i\pi (\lambda+\delta)} |\begin{array}{c} (2\beta, 2\delta),\;{}_1(a_j, e_j)_p\\ (2\beta, 2\delta), \;{}_1(b_j, f_j)_q\end{array} \right] 
\end{align}
where $H_{\;p,\; q}^{\;m,\; n}[z]$ is the well known H-function[1].\\
\textbf{Proof : } In order to establish the identity (2.1), we proceed as follows. \\
Denoting the left-hand of H-function by I, expressing the H-function with the help of its definition we have, 
\begin{equation}
I = 2\pi i \; \frac{1}{2\pi i}\;\int_L \theta(s)\;z^s\;  \frac{\Gamma(\beta-\delta s)\;\Gamma(1-\beta+\delta s)}{\Gamma(\alpha-\lambda s)\; \Gamma(1-\alpha+\lambda s)} \; ds
\end{equation}
where $\theta(s)$ is given by 
\begin{equation}
\theta(s) = \frac{\prod_{j=1}^m \Gamma(b_j-f_j s) \;\prod_{j=1}^n \Gamma(1-a_j+e_j s)}{\prod_{j=m+1}^{q} \Gamma(1-b_j+f_j s)\; \prod_{j=n+1}^p \Gamma(a_j-e_j s)}
\end{equation}
Using the results
\begin{equation}
\Gamma(\beta-\delta s) \; \Gamma(1-\beta+\delta s) = 2\pi\; \frac{\Gamma(2\beta- 2\delta s)\; \Gamma(1-2\beta+ 2\delta s)}{\Gamma(\h + \beta-\delta s)\; \Gamma(\h -\beta+\delta s)}
\end{equation} 
\begin{equation}
sin \pi z = \frac{\pi}{\Gamma(z) \Gamma(1-z)} = \frac{e^{i \pi z} - e^{-i \pi z}}{2 i}
\end{equation}
and 
\begin{equation}
cos \pi z = \frac{\pi}{\Gamma\left(\h -z \right)\; \Gamma\left(\h +z\right)} = \frac{e^{i \pi z} + e^{-i \pi z}}{2 }
\end{equation}
and after some algebra, we have 
\begin{align}
I & = \frac{1}{2\pi i}\;\int_L \theta(s)\; z^s\; \Gamma(2\beta-2\delta s)\;\Gamma(1- 2\beta+2\delta s)  \nonumber \\
&  \qquad . \left(e^{i\pi(\alpha-\lambda s)}- e^{-i\pi(\alpha-\lambda s)} \right)\; \left(e^{i\pi(\beta-\delta s)}+ e^{-i\pi(\beta-\delta s)}\right) ds  \nonumber \\
& = \frac{1}{2\pi i}\; \int_L \theta(s)\; z^s\; \Gamma(2\beta-2\delta s)\;\Gamma(1- 2\beta+2\delta s)  \nonumber \\
& \qquad .\left\{e^{i\pi(\alpha+\beta-\lambda s-\delta s)} + e^{i\pi(\alpha-\beta-\lambda s+\delta s)} \right.  \nonumber \\
& \qquad  -  \left. e^{-i\pi(\alpha-\beta-\lambda s+\delta s)} - e^{-i\pi(\alpha+\beta-\lambda s-\delta s)}\right\} ds
\end{align}
Now, breaking in to four parts and after some simplification, using the definition of H-function, we easily arrive at the right-hand side of (2.1). \\
This completes the proof of the identity (2.1).
\section{Special Cases}
\begin{enumerate}[(a)]
\item In (2.1), if we take $\delta=0$, we get, after some simplification, the identity obtained earlier  by Rathie[2].
\item In (2.1), if we take $\lambda=0$, we get, after some simplification, the identity obtained very recently by Rathie et al.[3].
\section{Concluding Remarks}
If we use the result (2.4), we get the following  identity for the H-function.
\begin{align}
 & H_{\;p+2,\; q+2}^{\;m+1,\; n+1} \left[ z |\begin{array}{c} (\beta, \delta), \;{}_1(a_j, e_j)_p, (\alpha, \lambda)\\ (\beta, \delta), \;{}_1(b_j, f_j)_q, (\alpha, \lambda)\end{array} \right] \nonumber\\
& \quad   =  \;  H_{\;p+4,\; q+4}^{\;m+2,\; n+2} \left[ z | \begin{array}{c} (2\beta, 2\delta), \left(\h + \alpha, \lambda\right), \;{}_1(a_j, e_j)_p, (2\alpha, 2\lambda), \left(\h + \beta, \delta \right) \\ (2\beta, 2\delta), \left(\h + \alpha, \lambda\right), \;{}_1(b_j, f_j)_q, (2\alpha, 2\lambda), \left(\h + \beta, \delta\right)\end{array} \right]   
\end{align}
Further in this, if we take $\delta=0$ and $\lambda=0$, we respectively get 
\begin{align}
 H_{\;p+1,\; q+1}^{\;m,\; n}& \left[ z |\begin{array}{c} \;{}_1(a_j, e_j)_p, (\alpha, \lambda)\\  \;{}_1(b_j, f_j)_q, (\alpha, \lambda)\end{array} \right] \nonumber\\
&  = \frac{1}{2\pi}\;H_{\;p+2,\; q+2}^{\;m+1,\; n+1} \left[ z |\begin{array}{c}  \left(\h + \alpha, \lambda\right), \;{}_1(a_j, e_j)_p, (2\alpha, 2\lambda) \\ \left(\h + \alpha, \lambda\right), \;{}_1(b_j, f_j)_q, (2\alpha, 2\lambda), \end{array} \right]  
\end{align}
and 
\begin{align}
 H_{\;p+1,\; q+1}^{\;m,\; n}& \left[ z |\begin{array}{c} (\beta, \delta), \;{}_1(a_j, e_j)_p\\ (\beta, \delta), \;{}_1(b_j, f_j)_q \end{array} \right] \nonumber\\
&  = 2\pi \; H_{\; p+2,\; q+2}^{\;m+1,\; n+1} \left[ z |\begin{array}{c} (2\beta, 2\delta),  \;{}_1(a_j, e_j)_p, \left(\h + \beta, \delta\right) \\ (2\beta, 2\delta),  \;{}_1(b_j, f_j)_q, \left(\h + \beta, \delta\right)\end{array} \right]
  \end{align}
 \end{enumerate}
 
\end{document}